\documentclass[11pt]{amsart}
\usepackage{amsmath,amssymb,amscd}

\newtheorem{thm}{Theorem}
\newtheorem{prop}{Proposition}
\newtheorem{coro}{Corollary}
\newtheorem{lemma}{Lemma}

\theoremstyle{definition}

\newtheorem{definition}{Definition}
\newtheorem{remark}{Remark}

\newtheorem*{eg}{Example}
\newtheorem{ex}{Exercise}

\def\tu{\textup}

\def\hom{\tu{Hom}}

\def\spec{\tu{Spec\,}}

\def\limit{\underset{t \to 0}{\lim} \, }

\def\Res{\textup{Res}}

\input xy
\xyoption{all}

\begin{document}

\title{Affine embeddings of a reductive group}
\author{David Murphy}
\date{\today} 
\address{Department of Mathematics and Computer Science\\
         Hillsdale College\\
         33 East College Street\\
         Hillsdale, MI  49242 }
\email{dmurphy@hillsdale.edu}
\begin{abstract}
We classify affine varieties with an action of a connected, 
reductive algebraic group such that the group is isomorphic 
to an open orbit in the variety.  This is accomplished by
associating a set of one-parameter subgroups of the group to 
the variety, characterizing such sets, and proving that sets 
of this type correspond to affine embeddings of the group.
Applications of this classification to the existence of
morphisms are then given.
\end{abstract}
\maketitle

\section{Introduction}\label{S:intro}

A basic problem in algebraic geometry is the study of algebraic actions.  Even actions with a dense orbit are not well understood.  Let $G$ be an algebraic group.  A \emph{quasihomogeneous variety} is a normal $G$-variety possessing an open orbit isomorphic to $G/H$, for some closed subgroup $H$ of $G$.  Toric varieties are an important family of quasihomogeneous varieties, where $G$ is an algebraic torus and $H$ is its trivial subgroup.  A partial classification of quasihomogeneous varieties was obtained in the important paper of Luna and Vust~\cite{LV}.  Their classification seeks to generalize that of toric varieties, but it is only feasible in special cases.  In this paper, we solve the equivariant classification problem for one case not 
covered in~\cite{LV}, namely the classification of affine quasihomogeneous varieties in which $H$ is trivial.

A \emph{$G$-embedding} is a normal $G$-variety $X$ that contains an open orbit $\Omega$ isomorphic to $G$.  The closed subvariety $\partial X = X - \Omega$ is called the \emph{boundary} of $X$.  Since $\Omega$ is an open $G$-orbit, $\partial X$ is a $G$-stable divisor of $X$ unless $\Omega = X$.The irreducible components of $\partial X$ are $G$-stable prime divisors of $X$.  Any toric variety is a $T$-embedding, in our terminology.  Similarly, the wonderful compactification of an adjoint group $G$ defined by De Concini and Procesi~\cite{CP1} is a $G$-embedding.  All of these examples have both a left and a right $G$-action.  Yet our definition of a $G$-embedding allows 
the consideration of $G$-varieties with only a left action of the group.  Hence, our definition of $G$-embeddings includes all of the biequivariant compactifications already in the literature~\cite{BCP},~\cite{B2},~\cite{CP1},~\cite{TE},~\cite{Tim2} and many more.  In this paper, we study affine $G$-embeddings and relate our results with those of Brion~\cite{B2} in the case of a biequivariant affine $G$-embedding in Section~\ref{SS:biequivariance}.

Suppose $X$ is an affine $G$-embedding and $x_0 \in X$ is a closed point in the open orbit $\Omega$.  In Section~\ref{S:onePS}, we define a set $\Gamma(X,x_0)$ of one-parameter subgroups of $G$ associated to the embedding $X$ and base point $x_0$.  Properties of such sets are collected in
Proposition~\ref{P:4properties}, which are then used in Theorem~\ref{T:A_X} to prove that $X$, as a $G$-embedding with basepoint $x_0$, is determined by its set $\Gamma(X,x_0)$.  Therefore, we
turn our attention to the classification of such sets in Sections~\ref{SS:1skeleton} through~\ref{SS:classification}.  We prove that sets $\Gamma(X,x_0)$ arising from affine $G$-embeddings are strongly convex lattice cones, in the sense of Definition~\ref{D:scrpc}, and that any strongly convex lattice cone determines an affine $G$-embedding in Theorem~\ref{T:classification}.  This generalizes the classification of affine toric varieties by strongly convex rational polyhedral cones in~\cite{TE}.  

Lastly, we explore the functoriality of our classification in Sections~\ref{SS:morphisms} and~\ref{SS:biequivariance}.  Specifically, Proposition~\ref{P:morphisms} states that equivariant morphisms of affine $G$-embeddings correspond to the inclusion of their associated cones, up to conjugation, analogous to the result for toric varieties.  Moreover, our classification reveals when an affine $G$-embedding $X$ has not only a left, but also a right $G$-action compatible with the identification of $G$ with the open orbit in $X$ (Proposition~\ref{P:biequivariance}).  Using this, we define a biequivariant resolution (Definition~\ref{D:biresolution}) of an arbitrary affine $G$-embedding and prove its universal property in Theorem~\ref{T:universalbi} of Section~\ref{SS:biequivariance}.

\subsection*{Notation}

We will always work over a ground field $k$, which we assume to be algebraically closed and of characteristic zero.  All algebraic groups are assumed to be linear and defined over $k$, and will be denoted by letters such as $G$ and $H$.  In particular, $G$ will refer to a connected, reductive algebraic group defined over $k$.  The symbol $T$ will always denote an algebraic torus, whether abstract or as a subgroup of $G$.  For an algebraic group $H$, $\mathfrak{X}_*(H)$ will
denote its set of one-parameter subgroups and $\mathfrak{X}^*(H)$ will denote the group of characters of $H$.

\subsection*{Acknowledgments}
I wish to express my sincere thanks to my advisors, Robert Fossum and William Haboush, for their encouragement and guidance as I conducted this work.

\section{One-parameter subgroups}\label{S:onePS}

Our primary method for describing and classifying affine $G$-embeddings $X$ is to make use of one-parameter subgroups of $G$.  We are interested in the limits, when they exist, of the one-parameter subgroups of $G$ in $X$.  One-parameter subgroups and their limits have been employed in a number of applications, including the Hilbert-Mumford criterion of stability~\cite{GIT}, the construction of the spherical building of the group $G$~\cite{GIT} and the Bialynicki-Birula decomposition of a smooth projective $T$-variety~\cite{BB}.  For our purposes, we will show that an affine $G$-embedding $X$ is determined by the set of one-parameter subgroups $\gamma$ of $G$ such that $\limit \gamma(t) x_0$ exists in $X$.

A \emph{one-parameter subgroup} of $G$ is a homomorphism of algebraic groups $\gamma : \mathbb{G}_m \to G$, which thus corresponds to a map $\gamma^\circ : k[G] \to k[t,t^{-1}]$.  Let $\mathfrak{X}_*(G)$ denote the set of one-parameter subgroups of $G$.  The group $G$ acts on $\mathfrak{X}_*(G)$ by conjugation, $g \bullet \gamma : t \mapsto g \gamma(t) g^{-1}$.  We will
denote the trivial one-parameter subgroup $t \mapsto e$ by $\varepsilon$.  Each one-parameter subgroup $\gamma \in \mathfrak{X}_*(G)$ determines a subgroup 
\begin{equation}\label{E:P(gamma)}
P(\gamma) = \{ g \in G : \gamma(t) g \gamma(t^{-1}) \in G_{k[[t]]} \}
\end{equation}
of $G$, which is parabolic if $G$ is reductive~\cite{GIT}.  In fact, every parabolic subgroup of a reductive group $G$ is of the form $P(\gamma)$ for some one-parameter subgroup $\gamma$ of $G$~\cite{Springer}.  We define an equivalence relation on the set of non-trivial one-parameter subgroups of $G$ by  
\begin{equation}\label{E:sim1ps}
\gamma_1 \sim \gamma_2 \text{ if and only if } \gamma_2(t^{n_2}) = g \gamma_1(t^{n_1}) g^{-1}
\end{equation}
for positive integers $n_1,n_2$ and an element $g \in P(\gamma_1)$, for all $t \in k^\times$.  Then the quotient $(\mathfrak{X}_*(G) - \{ \varepsilon \})/\sim$ is isomorphic to the spherical building of $G$~\cite{GIT},~\cite{Tits2}.

Every parabolic subgroup $P$ of $G$ defines a subset $\Delta_P(G) = \{ \gamma \in \mathfrak{X}_*(G) : P(\gamma) \supseteq P \}$ of $\mathfrak{X}_*(G)$.  Clearly $\gamma \in \Delta_{P(\gamma)}(G)$ for all $\gamma \in \mathfrak{X}_*(G)$, so $\mathfrak{X}_*(G) = \bigcup \Delta_P(G)$ where the union is over all parabolic subgroups of $G$.  In the spherical building of $G$, the images of the sets $\Delta_P(G)$ are simplices and constitute a ``triangulation'' of the building~\cite{GIT}.

The inclusion of $k[t,t^{-1}]$ in $k((t))$ allows us to view $\mathfrak{X}_*(G)$ as a subset of $G_{k((t))} = \hom_k(k[G],k((t)))$, the set of $k((t))$-points of $G$.  Let $\langle\gamma\rangle \in G_{k((t))}$ denote the point corresponding to the one-parameter subgroup $\gamma$.  The group $G_{k((t))}$ contains the subgroup $G_{k[[t]]}$, which consists of all $k((t))$-points of $G$ that have a specialization in $G$ as $t \to 0$.  The group $G_{k((t))}$ is the disjoint union of the double cosets of $G_{k[[t]]}$, as described by the Iwahori decomposition:

\begin{thm}[\cite{Iwahori}, Cartan--Iwahori Decomposition]\label{T:Iwahori} 
Let $G$ be a reductive algebraic group over $k$.  Every double coset of $G_{k((t))}$ with respect to the subgroup $G_{k[[t]]}$ is represented by a point of the type $\langle \gamma \rangle$, for some one-parameter subgroup $\gamma$ of $G$.  That is,
\begin{equation}\label{E:Iwahori}
G_{k((t))} = \bigcup_{\gamma \in \mathfrak{X}_*(G)} 
                     G_{k[[t]]} \langle \gamma \rangle G_{k[[t]]}
\end{equation}
Furthermore, each double coset is represented by a unique dominant one-parameter subgroup.
\end{thm}

This decomposition enables us to replace $k((t))$-points of $G$ with one-parameter subgroups.

\section{Limits of one-parameter subgroups}\label{SS:limits}

Let $X$ be a $G$-variety.  Each point $x$ of $X$ determines a morphism $\psi_x : G \to X$ by $\psi_x(g) = g \cdot x$.  For a point $x_0 \in X$ and a one-parameter subgroup $\gamma$ of $G$, we say $\limit \gamma(t) x_0$ exists in $X$ if $\psi_{x_0} \circ \gamma : \mathbb{G}_m \to X$ extends to a morphism $\tilde{\gamma} : \mathbb{A}^1 \to X$.  In this case, $\limit \gamma(t) x_0$ is defined to be $\tilde{\gamma}(0)$.  That is, the composition of $\psi_{x_0}^\circ : k[X] \to k[G]$ with $\gamma^\circ : k[G] \to k[t,t^{-1}]$ factors through $k[t]$, and the limit, $\limit \gamma(t)x_0$, is the $k$-point of $X$ corresponding to the composite $k[X] \to k[t] \to k$
sending $t \to 0$.  This is described by the diagrams:
\[
\xymatrix{
\mathbb{G}_m
  \ar[r]^\subset
  \ar[d]_\gamma
&
\mathbb{A}^1
  \ar@{-->}[d]^{\tilde{\gamma}_0}
&&
k[X]
  \ar[r]^{\psi_{x_0}^\circ}
  \ar@{-->}[d]_{\tilde{\gamma}_0^\circ}
&
k[G]
  \ar[d]^{\gamma^\circ}
\\
G
  \ar[r]_{\psi_{x_0}}
&
X
&&
k[t]
  \ar[r]_\subset
&
k[t,t^{-1}].
}
\]
Similarly, if $\lambda$ is a $k((t))$-point of $G$, then $\limit \lambda(t) x_0$ exists in $X$ means $\lambda^\circ|_{k[X]} : k[X] \to k[[t]]$.

The following lemma is used frequently hereafter.

\begin{lemma}\label{L:product_limit}
Suppose $\lambda \in G_{k((t))}$ and $\alpha \in G_{k[[t]]}$, so that $\alpha$ has specialization $\alpha_0 \in G_k$.  Let $X$ be an affine $G$-embedding with base point $x_0$.  Then $\limit [ \lambda(t) x_0 ]$ exists in $X$ if and only if $\limit [ \alpha(t) \lambda(t) x_0 ]$ exists, in which case
\begin{equation}\label{E:product_limit}
\limit [ \alpha(t) \lambda(t) x_0 ] = 
\alpha_0 \cdot \limit [\lambda(t) x_0].
\end{equation}
\end{lemma}

The proof is straightforward.

\begin{remark}\label{R:bi-limits}
Suppose $X$ is a biequivariant $G$-embedding, so $G$ has both a left and a right action on $X$ and $G$ may be identified with an open subvariety $\Omega$ which is stable for both actions.  Then we could amplify Lemma~\ref{L:product_limit} as follows:
If $\alpha,\beta \in G_{k[[t]]}$ and $\lambda \in G_{k((t))}$, then 
$\limit \lambda(t) x_0 \in X$ if and only if $\limit [ \alpha(t)
\lambda(t) \beta(t) x_0 ] \in X$, in which case
\[
\limit [ \alpha(t) \lambda(t) \beta(t) x_0 ] = 
\alpha_0 \cdot [ \limit \lambda(t) x_0 ] \cdot \beta'_0,
\]
where $\beta(t) \cdot x_0 = x_0 \cdot \beta'(t)$ for some $\beta' \in G_{k[[t]]}$
and $\alpha_0,\beta'_0 \in G_k$ denote the specializations of $\alpha,\beta'$, 
respectively.  However, we must be careful, for $\limit [\alpha(t)\lambda(t)
\beta(t) x_0]$ does not have to equal $\limit [\alpha_0 \lambda(t) \beta_0
x_0 ]$, where $\alpha_0,\beta_0$ are the specializations of $\alpha,\beta$~\cite{Kempf}.
\end{remark}

\begin{thm}[\cite{Kempf}, Theorem 1.4]\label{T:one-ps}
Let $X$ be an affine $G$-variety.  Suppose that $Y$ is a closed $G$-stable
subvariety of $X$ and that $x_0 \in X$ is a closed point such that the
closure of the orbit $Gx_0$ intersects $Y$.  Then there is a one-parameter
subgroup $\gamma$ of $G$ such that $\limit \gamma(t) x_0 \in Y$.
\end{thm}

\section{One-parameter subgroups of an affine $G$-embedding}\label{SS:Gamma(X)}

\begin{definition}\label{D:Gamma(X)}
Given a $G$-variety $X$ and a point $x_0 \in X$, define 
\begin{equation}\label{E:Gamma(X)}
\Gamma(X,x_0) := 
\{ \gamma \in \mathfrak{X}_*(G) : \limit \gamma(t) x_0 
                                  \text{ exists in } X \}.
\end{equation}
\end{definition}

We will be interested in the structure of such sets of one-parameter subgroups when $X$ is an affine $G$-variety and the orbit of $x_0$ in $X$ is open and isomorphic to $G$.  We call such an $x_0 \in X$ a base point.  Before we proceed, we make some immediate observations about such sets.

\begin{prop}\label{P:4properties}
Let $G$ be a connected reductive group.  Suppose $X$ is an affine $G$-embedding and $x_0 \in X$ is a base point.
\begin{enumerate}
  \item  [a.] If $x'_0 = hx_0$, then $\Gamma(X,x'_0) = h \Gamma(X,x_0) h^{-1}$.
  \item  [b.] If $\gamma \in \Gamma(X,x_0)$ and $\gamma \neq \varepsilon$, 
              then $\gamma^{-1} \not\in \Gamma(X,x_0)$.
  \item  [c.] If $T$ is any torus of $G$, then $\overline{Tx_0} \cong \overline{T}_\sigma$,
              where $\sigma = \Gamma(X,x_0) \cap \mathfrak{X}_*(T)$ is a strongly convex
              rational polyhedral cone in $\mathfrak{X}_*(T)$~\cite{TE}.
  \item  [d.] If $\gamma \in \Gamma(X,x_0)$ and $p \in P(\gamma)$, then $p \gamma(t) p^{-1}
              \in \Gamma(X,x_0)$, and moreover 
              \begin{equation}\label{E:GammaPs}
              \Gamma(X,x_0) = \bigcup_{P \subset G} 
                              P \bullet (\Gamma(X,x_0) \cap \Delta_P(G))
              \end{equation}
              where the union is taken over all parabolic subgroups $P$ of $G$.
  \item  [e.] The image of $\Gamma(X,x_0)$ in the spherical building is convex 
              (\cite{GIT}, Definition 2.10).
\end{enumerate}
\end{prop}

\begin{proof}
First, $\Gamma(X,x_0)$ depends on the base point $x_0$ as follows.  Suppose
$x'_0 \in \Omega$, so $x'_0 = h \cdot x_0$ for some unique $h \in G$ (because
$G \to \Omega$ is an isomorphism).  Then $\Gamma(X,x'_0) = h\Gamma(X,x_0)h^{-1}$
for if $\gamma \in \Gamma(X,x_0)$ (that is, if $\limit \gamma(t) x_0
\in X$), then $\limit (h\gamma(t)h^{-1}) x'_0 = \limit
h \gamma(t) h^{-1} h x_0 = \limit h\gamma(t) x_0 = h \limit
\gamma(t) x_0$, which exists in $X$.  Therefore $h\Gamma(X,x_0)h^{-1} \subset
\Gamma(X,x'_0)$.  By symmetry, since $x_0 = h^{-1} x'_0$, $h^{-1} \Gamma(X,x'_0)
h \subset \Gamma(X,x_0)$, so $\Gamma(X,x'_0) \subset h \Gamma(X,x_0) h^{-1}$.
Hence,
\begin{equation}\label{E:Gamma(X,hx)}
\Gamma(X,h \cdot x_0) = h\, \Gamma(X,x_0)\, h^{-1}
\end{equation}
for any $h \in G$. 

Second, as $X$ is affine, if $\gamma \in \Gamma(X,x_0)$ and $\gamma$ is not the
trivial one-parameter subgroup $\varepsilon : t \mapsto e$, then $\gamma^{-1} 
\not\in \Gamma(X,x_0)$.  Otherwise, if both $\limit \gamma(t) x_0$ and 
$\limit \gamma^{-1}(t) x_0$ exist in $X$, then the composition $\psi_{x_0} 
\circ \gamma : \mathbb{G}_m \to X$ extends to a morphism $\tilde{\gamma} : 
\mathbb{P}^1 \to X$, which must therefore be constant, so $\gamma = \varepsilon$.

Third, if $T$ is any torus of $G$, then $\overline{T x_0} \cong T_\sigma$, where 
$\sigma \subset \mathfrak{X}_*(T)$ is the strongly convex lattice cone 
$\Gamma(X,x_0) \cap \mathfrak{X}_*(T)$ from toric geometry.

Now suppose $\gamma \in \Gamma(X,x_0)$ and $p \in P(\gamma)$.  Then 
$p \cdot \gamma \cdot p^{-1}$ also belongs to $\Gamma(X,x_0)$, for $\limit 
[ (p \gamma(t) p^{-1}) x_0 ] = \limit [ p ( \gamma(t) p^{-1} \gamma(t^{-1}) ) 
\gamma(t) x_0 ] = p [\limit \gamma(t) p^{-1} \gamma(t^{-1}) ] [ \limit 
\gamma(t) x_0 ]$ exists in $X$ by Lemma~\ref{L:product_limit} and the definition of 
$P(\gamma)$.  Therefore, it is clear that
$\Gamma(X,x_0) = \bigcup_{P \subset G} P \bullet (\Gamma(X,x_0) \cap \Delta_P(G))$,
where the union is taken over all parabolic subgroups $P$ of $G$.

Lastly, if $\delta_1,\delta_2 \in (\Gamma(X,x_0) - \{ \varepsilon \})/\sim$,
then there are one-parameter subgroups $\gamma_1,\gamma_2 \in \Gamma(X,x_0)$ such 
that $\delta_i = [\gamma_i]$ is the equivalence class of $\gamma_i$.  The one-parameter
subgroups determine parabolic subgroups $P(\gamma_1)$ and $P(\gamma_2)$, whose intersection
contains a maximal torus $T$ of $G$.  Then $\gamma_1$ and $\gamma_2$ are equivalent to
one-parameter subgroups $\gamma_1',\gamma_2' \in \mathfrak{X}_*(T)$ and $\delta_i =
[\gamma_i']$.  By part 4, $\gamma_1', \gamma_2' \in \Gamma(X,x_0)$ as well.  Then 
$\gamma_1',\gamma_2' \in \Gamma(X,x_0) \cap \mathfrak{X}_*(T)$, which is the strongly 
convex rational polyhedral cone associated to the toric variety $\overline{T} \subset X$
by part 3.  As strongly convex rational polyhedral cones are convex, the line in 
$\mathfrak{X}_*(T)$ joining $\gamma_1$ and $\gamma_2$ is contained in $\Gamma(X,x_0) 
\cap \mathfrak{X}_*(T)$, and hence the line in the spherical building joining $\delta_1$ 
and $\delta_2$ is contained in the image of $\Gamma(X,x_0)$, so this image is 
semi-convex.  It is convex by part 2, which implies no pair of antipodal points of
the building can belong to the image of $\Gamma(X,x_0)$.
\end{proof}

Each $\gamma \in \mathfrak{X}_*(G)$ may be viewed as a $k((t))$-point of $G$.
In~\cite{LV}, a $G$-stable valuation $v_\lambda$ is associated to every $\lambda
\in G_{k((t))}$ in the following way.  As $\lambda$ is a $k((t))$-point of $G$,
we obtain a dominant morphism
\[
\xymatrix{
G \times \spec k((t)) \ar[r]^>>>>>>{1 \times \lambda}
&
G \times G \ar[r]^>>>>>\mu
&
G.
}
\]
This morphism induces an injection of fields $i_\lambda : k(G) \to  
\text{Frac}(k(G) \otimes_k k((t))) \to k(G)((t))$.  Then $v_t \circ i_\lambda : k(G)^\times 
\to \mathbb{Z}$ is a valuation of $k(G)$, where $v_t : k(G)((t))^\times \to 
\mathbb{Z}$ is the standard valuation associated to the order of $t$.  We 
define $v_\lambda = \dfrac{1}{n_\lambda} (v_t \circ i_\lambda)$, where 
$n_\lambda \in \mathbb{Z}$ is the largest positive number such that 
$(v_t \circ i_\lambda)(k(G)^\times) \subset n_\lambda \mathbb{Z}$ (except 
when $\lambda = \varepsilon$, in which case $v_\varepsilon(f) = 0$ or $\infty$ 
as $f(e) \neq 0$ or $= 0$, respectively).  This is $G$-stable by left 
translations, i.e., $v_\lambda(s \cdot f) = v_\lambda(f)$ for all $s \in G$, 
since $i_\lambda$ is clearly equivariant and $k(G)[[t]]$ is obviously stable 
for left translations by $G$ in $k(G)((t))$.  We include some of the properties
of these valuations that are proven in~\cite{LV} in the following lemma.

\begin{lemma}[\cite{LV}]\label{L:v_gamma}
\begin{enumerate}
  \item  [a.] Let $\gamma$ be a one-parameter subgroup of $G$.  For each $f \in k(G)$,
         there is an open subset $U \subset G$, depending only on $f$, such that
         \begin{equation}\label{E:inf}
         v_\gamma(f) = \inf_{s \in U} v_t(f(s \cdot \gamma(t)))
         \end{equation}
  \item  [b.] Let $\gamma_1,\gamma_2$ be one-parameter subgroups of $G$.  Then
         $v_{\gamma_1} = v_{\gamma_2}$ if and only if $\gamma_1 \sim \gamma_2$.
\end{enumerate}
\end{lemma}

\begin{proof}
Part 1 is Lemma 4.11.1 in~\cite{LV}, where $U = \{ s \in G : f(s) \neq 0 \}$.  The
second part is the result of Propositions 3.3 and 5.4 in~\cite{LV}.
\end{proof}

The sets of one-parameter subgroups $\Gamma(X,x_0)$ described in Definition~\ref{D:Gamma(X)}
are significant for the following reasons.  The first result, which will serve as our 
foundation for the classification theorem in Section~\ref{SS:classification}, is a 
uniqueness theorem which shows that an affine $G$-embedding $X$ with base point $x_0$ 
is determined by the set $\Gamma(X,x_0)$.  The second demonstrates that the prime 
divisors on the boundary of an affine $G$-embedding correspond to equivalence classes 
of edges of the set $\Gamma(X,x_0)$.  

\begin{thm}\label{T:A_X}
Let $G$ be a connected reductive group.  If $X$ is an affine $G$-embedding with base 
point $x_0$, then $X \cong \spec A_{\Gamma(X,x_0)}$, where
\begin{equation}\label{E:A_X}
A_{\Gamma(X,x_0)} := \{ f \in k[G] : v_\gamma(f) \geq 0 
                                     \text{ for all } \gamma \in \Gamma(X,x_0) \}.
\end{equation}
\end{thm}

\begin{proof}
The base point $x_0$ defines a morphism $\psi_{x_0} : g \mapsto g \cdot x_0$ from $G$ 
to $X$.  As both $G$ and $X$ are affine, $\psi_{x_0}$ corresponds to a homomorphism 
$\psi_{x_0}^\circ : k[X] \to k[G]$, which is injective since the image of $G$ is open 
in $X$ and $X$ is irreducible.  The image of $\psi_{x_0}^\circ$ lies in the 
subalgebra $A_{\Gamma(X,x_0)}$ since every $\gamma \in \Gamma(X,x_0)$ extends to 
a morphism $\tilde{\gamma} : \mathbb{A}^1 \to X$ so that $f(g \cdot \tilde{\gamma}(0))$ 
exists, which implies that $v_\gamma(f) = \inf_{s \in G_f} v_t(f(s \cdot \gamma(t))) 
\geq 0$ for all $f \in k[X]$.  We claim that $k[X] \cong A_{\Gamma(X,x_0)}$.  It suffices 
to show that every $f \in A_{\Gamma(X,x_0)}$ extends to a regular function on $X$ to prove 
that $k[X] \to A_{\Gamma(X,x_0)}$ is surjective and hence that $k[X] \cong 
A_{\Gamma(X,x_0)}$.

Suppose not and assume that $f \in A_{\Gamma(X,x_0)}$ is not in the image of 
$k[X]$.  Then $f$ fails to extend to a regular function on $X$, but it is 
defined on $\Omega = Gx_0$ by $f(g \cdot x_0) := f(g)$.  Let $P$ be the divisor 
of poles of $f$ in $X$.  Then $P$ is closed, has codimension one in $X$, and 
$P \subseteq \partial X$.  Let $D$ be an irreducible component of $\partial X$ 
and hence a closed subvariety of $X$.  Since $\partial X$ is $G$-stable and 
$G$ is connected (and so is irreducible), $D$ is a $G$-stable prime divisor 
since $e_G \in G$ fixes the generic point of the irreducible subvariety $D$.  
Therefore, Theorem~\ref{T:one-ps} provides a one-parameter subgroup $\gamma_D$ of 
$G$ with $\limit \gamma_D(t) x_0 \in D$, so $\gamma_D \in \Gamma(X,x_0)$.  
Yet $f \in A_{\Gamma(X,x_0)}$ implies that $v_{\gamma_D}(f) \geq 0$, so $f$ 
must be defined on the orbit $G [\limit \gamma_D(t) x_0 ] \subset D$.  
Therefore, $P \cap D$ is a closed subset of $D$ not equal to $D$, since $P$ 
does not contain $G [\limit \gamma_D(t) x_0 ] \subset D$ as 
$v_{\gamma_D}(f) \geq 0$.  Thus the codimension of $P$ in $X$, which is equal 
to the minimum of the codimensions of the $P \cap D$ as $D$ ranges over the 
irreducible components of $\partial X$, is at least $2$.  This is a contradiction.  
Hence every $f \in A_{\Gamma(X,x_0)}$ extends to a regular function on $X$, so 
is in the image of $\psi_{x_0}^\circ$.  Therefore, $\psi_{x_0}^\circ : k[X] \to
A_{\Gamma(X,x_0)}$ is an isomorphism, so $X \cong \spec A_{\Gamma(X,x_0)}$ as claimed.
Furthermore, the selected base point $x_0 \in X$ corresponds to the maximal ideal
$\mathfrak{m}_{x_0} = (\psi_{x_0}^\circ)^{-1}(\mathfrak{m}_e \cap A_{\Gamma(X,x_0)})$
of $k[X]$ as $\psi_{x_0}^\circ$ identifies $f \in A_{\Gamma(X,x_0)}$ with the 
unique extension of the function $f(g \cdot x_0) := f(g)$ to $X$.
\end{proof}

\begin{coro}\label{C:right-translation}
Let $X$ be an affine $G$-embedding with base point $x_0$.  If $x'_0 = hx_0$ is another
base point, then 
\begin{equation}\label{E:right-translation}
A_{\Gamma(X,h \cdot x_0)} = r_h(A_{\Gamma(X,x_0)}),
\end{equation}
where $r_h$ denotes right translation by $h$ in $k[G]$, $r_h(f)(x) = f(xh)$.
\end{coro}

\begin{proof}
Suppose $X$ is an affine $G$-embedding.  We have shown in (\ref{E:Gamma(X,hx)}) 
that the set $\Gamma(X,x_0)$ is determined by $X$ only up to conjugation, as any 
other base point is of the form $h \cdot x_0$ for a unique element $h \in G$ and
$\Gamma(X,h \cdot x_0) = h \Gamma(X,x_0) h^{-1}$.  We now show that 
$A_{\Gamma(X,h \cdot x_0)} = r_h(A_{\Gamma(X,x_0)})$.
Suppose $f \in r_h(A_{\Gamma(X,x_0)})$.  Then $f = r_h(f')$ for some $f' \in 
A_{\Gamma(X,x_0)}$, which means that $v_\gamma(f') \geq 0$ for all $\gamma \in 
\Gamma(X,x_0)$.  If $\gamma' \in \Gamma(X,h \cdot x_0) = h \Gamma(X,x_0) h^{-1}$, 
write $\gamma' = h \bullet \gamma$ for $\gamma \in \Gamma(X,x_0)$.
Using formula~\ref{E:inf}, one can easily show that $v_{h \bullet \gamma}(f) 
= v_\gamma(f') \geq 0$.  Therefore, $f \in A_{\Gamma(X,h \cdot x_0)}$, so $r_h(A_{\Gamma(X,x_0)}) 
\subset A_{\Gamma(X,h \cdot x_0)}$.  Likewise, 
$r_{h^{-1}}(A_{\Gamma(X,h \cdot x_0)}) \subset A_{\Gamma(X,x_0)}$, 
so $r_h(A_{\Gamma(X,x_0)}) = A_{\Gamma(X,h \cdot x_0)}$ as claimed.
\end{proof}


By Theorem~\ref{T:A_X} and Corollary~\ref{C:right-translation}, the
classification of affine $G$-embeddings is equivalent to the
characterization of such subsets of $\mathfrak{X}_*(G)$ that are
obtained from affine $G$-embeddings.  In order to classify admissible
subsets $\Gamma$, we will explore the properties of the sets
$\Gamma(X,x_0)$ for arbitrary affine $G$-embeddings $X$ with choice of
base point $x_0$ in the following sections.

\section{The one-skeleton of $\Gamma(X,x_0)$}\label{SS:1skeleton}

If $\sigma$ is a strongly convex rational polyhedral cone in 
$\mathfrak{X}_*(T)_\mathbb{R}$, let $\sigma(1)$ denote the set of rays
of $\sigma$, $\sigma(1) = \{ \tau < \sigma : \dim \tau = 1 \}$.  This
is called the one-skeleton of the cone.  By Proposition~\ref{P:4properties}, 
for each maximal torus $T$ of $G$, $\Gamma(X,x_0) \cap \mathfrak{X}_*(T)$ 
is a strongly convex rational polyhedral cone in $\mathfrak{X}_*(T)$.

\begin{definition}\label{D:one-skeleton}
Let $X$ be an affine $G$-embedding with base point $x_0$.
The \emph{one-skeleton} of the set $\Gamma(X,x_0)$ is
\begin{equation}\label{E:Gamma(1)}
\Gamma_1(X,x_0) := \bigcup_T \, [\Gamma(X,x_0) \cap \mathfrak{X}_*(T)](1),
\end{equation}
which is the set of extremal rays of $\Gamma(X,x_0)$.
\end{definition}

Recall the equivalence relation on one-parameter subgroups defined in Equation~\ref{E:sim1ps}.  Equivalence classes under this relation are now given a geometric interpretation.

\begin{prop}\label{P:Gamma(1)}
There is a bijection between $\Gamma_1(X,x_0)/\sim$ and the finite set
of prime divisors of $X$ contained in $\partial X$.
\end{prop}

\begin{proof}
Write $X = \Omega \cup D_1 \cup D_2 \cup \cdots \cup D_r$, where
$D_1,\dots,D_r$ are the irreducible components of $\partial X$.  Thus
each $D_i$ is a $G$-stable prime divisors of $X$.

Suppose that $\rho \in \Gamma_1(X,x_0)$.  Then there is a maximal torus $T$ 
of $G$ such that $\rho \in [\Gamma(X,x_0) \cap \mathfrak{X}_*(T)](1)$ is a 
ray of the strongly convex rational polyhedral cone $\Gamma(X,x_0) \cap 
\mathfrak{X}_*(T)$ corresponding to $\overline{T} \subset X$.  From the 
description of $T$-stable divisors in toric varieties in~\cite{Fulton}, the 
ray $\rho$ corresponds to a prime divisor $D_\rho$ in $\overline{T}$.  Thus 
$D_\rho$ is an irreducible $T$-stable subvariety of $\overline{T}$, so 
$\overline{G \cdot D_\rho}$ is an irreducible $G$-stable subvariety of $X$, 
as both $G$ and $D_\rho$ are irreducible.  Moreover, $\overline{G \cdot D_\rho}$ 
is contained in $\partial X$, so there is a $G$-stable prime divisor $D_i$ of $X$ 
such that $\overline{G \cdot D_\rho} \subseteq D_i \subset X$.  However, 
$\text{codim}_X(\overline{G \cdot D_\rho}) = \dim X - \dim
\overline{G \cdot D_\rho} = \dim G - ( \dim G + \dim D_\rho - \dim T ) = 
\dim T - \dim D_\rho = \text{codim}_{\overline{T}}(D_\rho) = 1$.  Thus, as $\overline{G
\cdot D_\rho}$ is irreducible and of codimension one in $X$, and $\overline{G \cdot 
D_\rho} \subset D_i \subsetneq X$, where $D_i$ is also irreducible and of codimension 
one, we conclude that $\overline{G \cdot D_\rho} = D_i$.  Therefore, there is a map
$\varphi : \Gamma_1(X,x_0) \to \{ D_1,D_2,\dots,D_r \}$ given by
$\varphi(\rho) = \overline{G \cdot D_\rho}$.

The map $\varphi$ is surjective, as any of the prime divisors $D_i$ of $X$ are 
closed $G$-subvarieties, and thus contain the limit point of some one-parameter 
subgroup $\gamma \in \Gamma(X,x_0)$ by Theorem~\ref{T:one-ps}.  This
$\gamma$ is a one-parameter subgroup of some maximal torus $T$ of $G$, so 
$\overline{T} \cap D_i \neq \emptyset$ is a $T$-stable divisor of $\overline{T}$.  
Hence, there is a prime divisor $D_\rho$ of $\overline{T}$ corresponding 
to a ray $\rho \in \Gamma(X,x_0) \cap \mathfrak{X}_*(T)$ such that $D_\rho \subset 
\overline{T} \cap D_i$.  Thus $\overline{G \cdot D_\rho} \subset D_i$,
from which we conclude $\overline{G \cdot D_\rho} = D_i$ as before.

Furthermore, $\varphi$ respects the equivalence relation $\sim$ described in 
(\ref{E:sim1ps}), for if $\gamma_\rho$ denotes the first lattice point in 
$\mathfrak{X}_*(T)$ along $\rho$ and $\gamma_{\rho_1} \sim \gamma_{\rho_2}$, then 
$\limit \gamma_{\rho_1}(t) x_0$ and $\limit \gamma_{\rho_2}(t) x_0$
belong to the same $G$-orbit in $X$, and hence to the same prime divisor $D$ of $X$,
as $D$ is $G$-stable.  Therefore, $\overline{G \cdot D_{\rho_1}} = D = \overline{G
\cdot D_{\rho_2}}$, from the discussion above.  Hence $\varphi$ induces a 
well-defined map $\widetilde{\varphi} : \Gamma_1(X,x_0)/\sim\, \to 
\{ D_1,D_2,\dots,D_r \}$, which is surjective.

We prove that $\widetilde{\varphi}$ is an injection.  Let $\rho \in \Gamma_1(X,x_0)$ and 
let $\gamma_\rho$ denote the first lattice point in $\mathfrak{X}_*(G)$ along 
$\rho$ as above.  The ideal $\Gamma(X,\mathcal{O}(\overline{G \cdot D_\rho})) 
= \{ f \in k[X] : f = 0 \text{ on } \overline{G \cdot D_\rho} \}$ is
equal to $k[X] \cap \mathfrak{m}_{v_{\gamma_\rho}}$, where
$\mathfrak{m}_{v_{\gamma_\rho}} = \{ f \in k(G) : v_{\gamma_\rho}(f) > 0 \}$ 
is the maximal ideal of the valuation ring $\mathcal{O}_{v_{\gamma_\rho}} = 
\{ f \in k(G) : v_{\gamma_\rho}(f) \geq 0 \}$.  For $D_\rho = 
\overline{T \cdot z_\rho}$, where $z_\rho = \limit \gamma_\rho(t) x_0 \in 
\overline{T}$, so that $\overline{G \cdot D_\rho} = \overline{G \cdot z_\rho}$ 
as well.  Thus, if $f \in k[X] \cap \mathfrak{m}_{v_{\gamma_\rho}}$, then 
$f(z_\rho) = f(\limit \gamma_\rho(t) x_0) = \limit t^{v_{\gamma_\rho}(f)} u = 0$, 
where $u \in k[X] \cap \mathcal{O}_{v_{\gamma_\rho}}^\times$, since 
$v_{\gamma_\rho}(f) > 0$.  Therefore $k[X] \cap \mathfrak{m}_{v_{\gamma_\rho}} 
\subset \Gamma(X,\mathcal{O}(\overline{G \cdot D_\rho}))$, where both are
prime ideals in $k[X]$ and the latter is of height one.  Hence they
are equal.  Therefore, $v_{\gamma_\rho}$ is the valuation of the prime
divisor $\overline{G \cdot D_\rho}$ in $X$.  Suppose $\rho_1,\rho_2
\in \Gamma_1(X,x_0)$ such that $\gamma_{\rho_1} \not\sim \gamma_{\rho_2}$.  
Then $v_{\gamma_{\rho_1}} \neq v_{\gamma_{\rho_2}}$, so 
$\overline{G \cdot D_{\rho_1}} \neq \overline{G \cdot D_{\rho_2}}$.
Hence $\Gamma_1(X,x_0)/\sim\, \to \{D_1,D_2,\dots,D_r\}$ is injective.
Therefore, $\Gamma_1(X,x_0)/\sim\, \to \{ D_1,D_2,\dots,D_r \}$ is a bijection.
\end{proof}

\section{Kempf states}\label{SS:Gstates}

Our sets $\Gamma(X,x_0)$ of one-parameter subgroups associated to an affine
$G$-embedding are identical to sets arising in geometric invariant theory as 
studied by Mumford~\cite{GIT}, Kempf~\cite{Kempf} and Rousseau~\cite{Rousseau}.
In~\cite{Kempf}, Kempf describes such sets in terms of bounded, admissible 
states as follows.

\begin{definition}\label{D:Kempfstate}
A \emph{state} $\Xi$ is an assignment of a nonempty subset $\Xi(R) \subset 
\mathfrak{X}^*(R)$ to each torus $R$ of $G$ so that the image of $\Xi(R_2)$ in 
$\mathfrak{X}^*(R_1)$ under the restriction map $\mathfrak{X}^*(R_2) \to 
\mathfrak{X}^*(R_1)$ is equal to $\Xi(R_1)$ whenever $R_1 \subset R_2$ are tori of $G$.

For each $k$-point $g$ of $G$, we have maps $g_! : \mathfrak{X}^*(g^{-1}Rg) 
\to \mathfrak{X}^*(R)$ defined by $(g_! \chi)(r) = \chi(g^{-1}rg)$ for each
torus $R$ of $G$.  We define the {\it conjugate state} $g_*\Xi$ by the
formula $(g_*\Xi)(R) = g_!\Xi(g^{-1}Rg)$ for each torus $R$ of $G$.  With
this notion in mind, we say that a state $\Xi$ is \emph{bounded} if, for
each torus $R$ of $G$, $\bigcup_{g \in G_k} g_*\Xi(R)$ is a finite set of
characters of $R$.

Finally, any state defines a function $\mu(\Xi)$ on $\mathfrak{X}_*(G)$ by
\begin{equation}\label{E:mu(Xi)}
\mu(\Xi,\gamma) = 
\min_{\chi \in \Xi(\gamma(\mathbb{G}_m))} \langle \chi,\gamma \rangle
\end{equation}
called the \emph{numerical function} of $\Xi$.  We say $\Xi$ is 
\emph{admissible} if its numerical function satisfies $\mu(\Xi,\gamma) =
\mu(\Xi,p \bullet \gamma)$ for all $p \in P(\gamma)$, where $P(\gamma)$ is
the parabolic subgroup of $G$ associated to $\gamma$.  We will refer to 
bounded admissible states as \emph{Kempf states}.
\end{definition}

\begin{remark}\label{R:restrictedstate}
Suppose $H$ is a closed subgroup of a group $G$ and that $\Xi$ is a Kempf state
for $G$.  Then $\Res^G_H \Xi$, which assigns to any torus $R$ of $H$ the set of
characters $\Xi(R)$ (for $R$ is also a torus of $G$), is a Kempf state for $H$,
as all of the compatibility conditions are clearly inherited from $G$.
\end{remark}

With these terms defined, we return to the situation of an affine $G$-scheme 
$X$.  For each closed $G$-subscheme $Y$ of $X$, define $\Gamma_Y(X,x_0)$ to 
be the set $\{ \gamma \in \Gamma(X,x_0) : \limit \gamma(t) \cdot x_0 
\text{ exists in } Y \}$.


\begin{thm}[\cite{Kempf}, Lemma 3.3]\label{T:gammastate}
Let $x_0$ be a $k$-point of an affine $G$-variety $X$.  Let $Y$ be a closed 
$G$-subvariety of $X$ not containing $x_0$.  Then there are Kempf states 
$\Xi_{X,x_0}$ and $\Upsilon^Y_{X,x_0}$ such that 
\begin{enumerate}
  \item  [a.] $\Gamma(X,x_0) = \{ \gamma \in \mathfrak{X}_*(G) : 
                             \mu(\Xi_{X,x_0},\gamma) \geq 0 \}$,
  \item  [b.] $\Gamma_Y(X,x_0) = \{ \gamma \in \Gamma(X,x_0) : 
                               \mu(\Upsilon^Y_{X,x_0},\gamma) > 0 \}$.
\end{enumerate}
\end{thm}

While not including the proof, we indicate the construction of the
Kempf state $\Xi_{X,x_0}$ given in~\cite{Kempf}.  By the embedding
theorem (Lemma 1.1 in~\cite{Kempf}), there is a $G$-representation $V$
and an equivariant closed embedding $X \hookrightarrow V$.  Identify $X$
with its image in $V$.  Since $\psi_{x_0} : G \to X$ is an isomorphism
onto the open orbit $\Omega \subset X$, we may ensure $x_0$ is not
zero in $V$.  As $X$ is a closed $G$-subvariety of $V$, $\Gamma(X,x_0)
= \Gamma(V,x_0)$, so we may assume $X = V$ is a $G$-representation.  
We define the state $\Xi_{V,x_0}$ of $x_0$ in the representation $V$
as follows.  Let $R$ be a torus of $G$.  Let $V = \bigoplus_{\chi \in
\mathfrak{X}^*(R)} V_\chi$ be the eigendecomposition of $V$ with
respect to the torus $R$ and let $\text{proj}_{V_\chi}(x_0)$ be the
projection of $x_0$ on the weight space $V_\chi$.  Set
\[
\Xi_{V,x_0}(R) = 
\{ \chi \in \mathfrak{X}^*(R) : \text{proj}_{V_\chi}(x_0) \neq 0 \},
\]
for each torus $R$ in $G$.  Then $\Xi_{V,x_0}$ is the Kempf state
associated to the set $\Gamma(X,x_0) = \Gamma(V,x_0)$.

Given a Kempf state $\Xi$, define the set $\Xi^\vee \subset \mathfrak{X}_*(G)$ by
\begin{equation}\label{E:XiV}
\Xi^\vee := \{ \gamma \in \mathfrak{X}_*(G) : \mu(\Xi,\gamma) \geq 0 \}.
\end{equation}
By Theorem~\ref{T:gammastate}, every collection of one-parameter subgroups $\Gamma(X,x_0)$ arising from an affine $G$-embedding is of the form $\Xi^\vee$ for some Kempf state $\Xi$.  Furthermore, observe that $\sum g \cdot v_\chi$ is the eigendecomposition of $g \cdot x_0$ with respect to $R$ whenever $\sum v_\chi$ is the eigendecomposition of $x_0$ with respect to $T = g^{-1}Rg$, so
\begin{equation}\label{E:g*Xi}
\Xi_{V,g \cdot x_0}(R) = g_!(\Xi_{V,x_0}(g^{-1}Rg)) = g_*\Xi_{V,x_0}(R)
\end{equation}
for all $g \in G_k$.  Thus, if $\Gamma(X,x_0) = \Xi^\vee$, then $g \bullet \Gamma(X,x_0) = \Gamma(X,g \cdot x_0) = (g_*\Xi)^\vee$.

Therefore, given an affine $G$-embedding $X$ and a choice of base point $x_0 \in \Omega$, we obtain a Kempf state $\Xi_{X,x_0}$ such that $\Gamma(X,x_0) = \{ \gamma \in \mathfrak{X}_*(G) : \mu(\Xi_{X,x_0},\gamma) \geq 0 \}$.  For each closed $G$-subvariety $Y$ of $X$, we obtain another Kempf state $\Upsilon^Y_{X,x_0}$ such that $\Gamma_Y(X,x_0) = \{ \gamma \in \Gamma(X,x_0) : 
\mu(\Upsilon^Y_{X,x_0},\gamma) > 0 \}$.

\section{Strongly convex lattice cones}\label{SS:cones}

One-parameter subgroups $\gamma_1, \gamma_2 \in \mathfrak{X}_*(G)$ are equivalent, as defined in Equation~\ref{E:sim1ps}, if and only if, for all $t \in k^\times$,
\[
\gamma_2(t^{n_2}) = g \gamma_1(t^{n_1}) g^{-1}
\]
for some positive integers $n_1,n_2$ and some element $g \in P(\gamma_1)$.  In this case we write $\gamma_1 \sim \gamma_2$.

\begin{definition}\label{D:scrpc}
We say a subset $\Gamma \subset \mathfrak{X}_*(G)$ is \emph{saturated} (with respect to the equivalence relation (\ref{E:sim1ps}) of one-parameter subgroups) if, whenever $\gamma_1 \sim \gamma_2$ and $\gamma_1 \in \Gamma$, then $\gamma_2 \in \Gamma$.  $\Gamma$ is called
a \emph{lattice cone} of one-parameter subgroups of $G$ if $\Gamma$ is saturated and the quotient $\Gamma_1/\sim$ of the one-skeleton of $\Gamma$ (\ref{E:Gamma(1)}) is a finite set.  A lattice cone $\Gamma$ is called a \emph{convex lattice cone} if there is a Kempf state $\Xi$ such that $\Gamma = \Xi^\vee$.  Additionally, $\Gamma$ is a \emph{strongly convex lattice cone} if it is a convex lattice cone and $\gamma,\gamma^{-1} \in \Gamma$ if and only if $\gamma$ is the trivial
one-parameter subgroup of $G$.
\end{definition}

The strong convexity condition implies that the elements of $\Xi(R)$ generate $\mathfrak{X}^*(R)$ as a group for all tori $R$ in $G$.

Our terminology is compatible with that of toric geometry.  If $\Gamma$ is a strongly convex lattice cone as above, then each $\Gamma \cap \mathfrak{X}_*(T)$ is one in the sense of toric geometry.  For if $\Gamma$ is a strongly convex lattice cone, $\Gamma = \{ \gamma \in \mathfrak{X}_*(G) : \mu(\Xi,\gamma) \geq 0 \}$ for some Kempf state $\Xi$.  Thus each $\Gamma \cap \mathfrak{X}_*(T)$ is the intersection of finitely many half-spaces, $\Gamma \cap \mathfrak{X}_*(T) = \bigcap_{\chi \in \Xi(T)} \{ v \in \mathfrak{X}_*(T) : \langle \chi,v \rangle \geq 0 \}$, so it is a convex lattice cone in $\mathfrak{X}_*(T)$.  Strong convexity follows as the same condition is required of $\Gamma$.

\begin{lemma}\label{L:Gammais}
If $X$ is an affine $G$-embedding with base point $x_0 \in \Omega$, then the set $\Gamma(X,x_0) = \{ \gamma \in \mathfrak{X}_*(G) : \limit \gamma(t) x_0 \text{ exists in } X \}$ is a strongly convex lattice cone.
\end{lemma}

\begin{proof}
Suppose $X$ is an affine $G$-embedding and $x_0 \in X$ is a base point.  Then $\Gamma(X,x_0) = \{ \gamma \in \mathfrak{X}_*(G) : \limit \gamma(t)x_0 \text{ exists in } X \}$.   By Theorem~\ref{T:gammastate}, $\Gamma(X,x_0) = \Xi^\vee$ for some Kempf state.  This, together with Proposition~\ref{P:4properties}.d, implies that $\Gamma(X,x_0)$ is saturated. Moreover, $\Gamma_1(X,x_0)/\sim$ is finite by Proposition~\ref{P:Gamma(1)}, so $\Gamma(X,x_0)$ is a convex lattice cone.  Therefore, it is a strongly convex lattice cone by Proposition~\ref{P:4properties}.b.
\end{proof}

\section{The classification of affine $G$-embeddings}\label{SS:classification}

Suppose $\gamma \in \mathfrak{X}_*(G)$.  Recall that $v_\gamma$ is the valuation $\dfrac{1}{n_\gamma}(v_t \circ i_\gamma)$ of $k(G)$, where $i_\gamma : k(G) \to k(G)((t))$ is the injection of fields filling the commutative diagram
\[
\xymatrix{
k[G] \ar[r]^>>>>>{\mu^\circ} \ar[d]^\subset
&
k[G] \otimes k[G] \ar[rr]^{id_{k[G]} \otimes \gamma^\circ}
&&
k[G] \otimes k((t)) \ar[d]^\subset
\\
k(G) \ar[rrr]_{i_\gamma}
&&&
k(G)((t))
}
\]
$v_t : k(G)((t))^\times \to \mathbb{Z}$ is the standard valuation, and $n_\gamma$ is the largest positive integer such that $(v_t \circ i_\gamma)(k(G)^\times) \subset n_\gamma\mathbb{Z}$.  Let $\mathcal{O}_{v_\gamma} = \{ f \in k(G) : v_\gamma(f) \geq 0 \}$ be the valuation ring in $k(G)$ associated to $v_\gamma$.

\begin{lemma}\label{L:A(W)}
Let $T$ be a maximal torus of $G$.  If $\gamma_1,\gamma_2 \in \mathfrak{X}_*(T) \subset \mathfrak{X}_*(G)$, then
\begin{equation}\label{E:O1capO2inO12}
(k[G] \cap \mathcal{O}_{v_{\gamma_1}}) \cap (k[G] \cap \mathcal{O}_{v_{\gamma_2}}) 
\subseteq k[G] \cap \mathcal{O}_{v_{\gamma_1+\gamma_2}}.
\end{equation}
\end{lemma}

\begin{proof}
Suppose $\gamma_1,\gamma_2$ are one-parameter subgroups of a maximal torus $T$ of $G$.  Then $\gamma_1 + \gamma_2 \in \mathfrak{X}_*(T)$ is also a one-parameter subgroup of $G$ contained in $T$.  The valuations $v_{\gamma_1},v_{\gamma_2}$ and $v_{\gamma_1+\gamma_2}$ are obtained from the 
homomorphisms $i_{\gamma_1},i_{\gamma_2}$ and $i_{\gamma_1+\gamma_2}$, so it is enough to prove that $i_{\gamma_1+\gamma_2}(k[G] \cap \mathcal{O}_{v_{\gamma_1}} \cap \mathcal{O}_{v_{\gamma_2}}) \subset k(G)[[t]]$ to prove the lemma.  Yet 
\[
\xymatrix{
k[G] \cap \mathcal{O}_{v_{\gamma_i}} \ar[rrr] \ar[d]^\subset
&&&
k[G] \otimes k[[t]] \ar[d]^\subset
\\
k[G] \ar[r]_{\mu^\circ}
&
k[G] \otimes k[G] \ar[rr]_{id_{k[G]} \otimes \gamma_i^\circ}
&&
k[G] \otimes k((t))
}
\]
for $i=1,2$ and
\[
\xymatrix{
k[G] \ar[r]^{\mu^\circ}
&
k[G] \otimes k[G] \ar[rr]^{id_{k[G]} \otimes (\gamma_1+\gamma_2)^\circ}
                  \ar[d]_{id_{k[G]} \otimes \mu^\circ}
&&
k[G] \otimes k((t))
\\
&
k[G] \otimes k[G] \otimes k[G] \ar[rr]_{id_{k[G]} \otimes \gamma_1^\circ \otimes \gamma_2^\circ}
&&
k[G] \otimes k((t)) \otimes k((t)) \ar[u]_{id_{k[G]} \otimes m_{23}}
}
\]
implies that $i_{\gamma_1 + \gamma_2}(k[G] \cap \mathcal{O}_{v_{\gamma_1}} \cap \mathcal{O}_{v_{\gamma_2}}) \subset k[G] \otimes k[[t]]$.  Thus, $(k[G] \cap \mathcal{O}_{v_{\gamma_1}}) \cap (k[G] \cap \mathcal{O}_{v_{\gamma_2}}) \subseteq k[G] \cap \mathcal{O}_{v_{\gamma_1+\gamma_2}}$ as claimed. 
\end{proof}

\begin{thm}\label{T:classification}
Affine $G$-embeddings are classified by strongly convex lattice cones of one-parameter subgroups in $\mathfrak{X}_*(G)$, up to conjugation.  Conjugation corresponds to the change of base point in the embedding.
\end{thm}

\begin{proof}
Assume $X$ is an affine $G$-embedding.  The selection of a base point $x_0$ in $X$ 
uniquely determines a set $\Gamma(X,x_0) \subset \mathfrak{X}_*(G)$, which is a strongly 
convex lattice cone by Lemma~\ref{L:Gammais}.  By Theorem~\ref{T:A_X}, the 
set $\Gamma(X,x_0)$ determines the affine $G$-embedding $X$ and the selected base 
point $x_0$ via an isomorphism $\psi_{x_0} : k[X] \cong A_{\Gamma(X,x_0)}$ such
that $\mathfrak{m}_{x_0} = \psi_{x_0}^{-1}(\mathfrak{m}_e \cap A_{\Gamma(X,x_0)})$.
By Corollary~\ref{C:right-translation} and formula (\ref{E:Gamma(X,hx)}), the selection
of a different base point $h \cdot x_0$ is equivalent to conjugating the cone 
$\Gamma(X,h \cdot x_0) = h \Gamma(X,x_0) h^{-1}$ and a corresponding right translation
of the algebra $A_{\Gamma(X,h \cdot x_0)} = r_h(A_{\Gamma(X,x_0)})$.  Thus each
affine $G$-embedding determines a strongly convex lattice cone, modulo 
conjugation in $\mathfrak{X}_*(G)$, which in turn recovers the variety up to isomorphism.

Conversely, we prove that every strongly convex lattice cone $\Gamma \subset 
\mathfrak{X}_*(G)$ determines an affine $G$-embedding $X_\Gamma$ such that 
$\Gamma(X_\Gamma,x_0) = \Gamma$ for some choice of base point $x_0 \in X_\Gamma$.  
Assume $\Gamma$ is a strongly convex lattice cone, so $\Gamma = \Xi^\vee$ for 
some Kempf state $\Xi$.  Define $A_\Gamma$ as in (\ref{E:A_X}) to be the 
subalgebra of $k[G]$ given by
\[
A_\Gamma = \{ f \in k[G] : v_\gamma(f) \geq 0, \text{ for all } \gamma
\in \Gamma \} = k[G] \cap \bigcap_{\gamma \in \Gamma} \mathcal{O}_{v_\gamma},
\]
where $\mathcal{O}_v$ is the valuation subring $\{ f \in k(G) : v(f) \geq 0 \}$ of 
$k(G)$.  Using this second description, since $k[G]$ and all of the valuation rings 
$\mathcal{O}_{v_\gamma}$ are integrally closed in $k(G)$, their intersection 
$A_\Gamma$ is an integrally closed domain.  Furthermore, $A_\Gamma$ is left-invariant, 
since $k[G]$ is and $v_\gamma(g \cdot f) = v_\gamma(f)$ for all $f \in k(G)$ and 
$g \in G$ implies that the valuation rings $\mathcal{O}_{v_\gamma}$ are $G$-stable 
as well.  Hence $A_\Gamma$ is an integrally closed left-invariant subalgebra of $k[G]$.  
It remains to prove that $A_\Gamma$ is finitely generated over $k$ and that the 
variety $X_\Gamma = \spec A_\Gamma$ contains $G$ as an open orbit.

Let $\gamma_1,\dots,\gamma_m$ be a set of representatives for the equivalence classes 
in $\Gamma_1/\sim$, which is a finite set as $\Gamma$ is a lattice cone.  By 
Proposition~\ref{P:Gamma(1)} and Lemma~\ref{L:A(W)}, 
\[
A_\Gamma 
= k[G] \cap \bigcap_{\gamma \in \Gamma} \mathcal{O}_{v_\gamma}
= \bigcap_{\gamma \in \Gamma} (k[G] \cap \mathcal{O}_{v_\gamma})
= \bigcap_{i=1}^m (k[G] \cap \mathcal{O}_{v_{\gamma_i}})
\]
since any $\gamma \in \Gamma$ belongs to some $\mathfrak{X}_*(T)$, in which case it is 
a sum of elements from $[\Gamma \cap \mathfrak{X}_*(T)](1)$.  Let $f_1,\dots,f_n$ be a 
finite set of generators for the algebra $k[G]$ over $k$, so $k[G] = k[f_1,\dots,f_n]$.  
By Gordan's Lemma, the set 
\[
C = \{ (c_\ell) \in \mathbb{Z}^n : c_1,\dots,c_n \geq 0, 
\sum_{j=1}^n c_j v_{\gamma_i}(f_j) \geq 0 \text{ for all } i = 1,\dots,m \}
\]
is finitely generated.  Thus the ``monomials'' $f_1^{c_1} f_2^{c_2} \cdots f_n^{c_n}$ 
in the generators $f_1,\dots,f_n$ of $k[G]$ associated to the generators $(c_1,\dots,c_n)$ 
of $C$ admit a finite set of generators of $A_\Gamma$, as any element of $A_\Gamma$ is 
an algebraic combination of such monomials.

Therefore $X_\Gamma = \spec A_\Gamma$ is a normal affine $G$-variety with an open 
$G$-orbit, since $A_\Gamma \subset k[G]$ implies $f : G \to X_\Gamma$ is dominant, so
$f(G)$ is an open orbit in $X_\Gamma$ (\cite{Springer}, Theorem 1.9.5).  It only 
remains to show that $f(G)$ is isomorphic to $G$ as an orbit in $X_\Gamma$.  Consider, 
for each torus $R$ of $G$, the image of $A_\Gamma$ under the homomorphism $k[G] \to k[R]$.  
By construction of $A_\Gamma = \{ f \in k[G] : v_\gamma(f) \geq 0 \text{ for all } \gamma 
\in \Gamma \}$ and the decomposition $A_\Gamma = \bigoplus_{\chi \in \Gamma^\vee(R)} 
A^R_\chi$, where $A^R_\chi = \{ f \in A : f(xr) = \chi(r)f(x) \text{ for all } r \in R \}$, 
it is evident that the image of $A_\Gamma$ in $k[R]$ is the $k$-monoid algebra 
$k[\Gamma^\vee(R)]$.  Hence we have the following commutative diagrams:
\[
\xymatrix{
A_\Gamma \ar[rr]^\subset \ar[d]_{\text{onto}}
&&
k[G] \ar[d]^{\text{onto}}
&
X_\Gamma
&&
G \ar[ll]_f
\\
k[\Gamma^\vee(R)] \ar[rr]_{\subset}
&&
k[R] = k[\mathfrak{X}^*(R)]
&
R_{\Gamma^\vee(R)} \ar[u]^{\text{closed}}
&&
R \ar[u]_{\text{closed}} \ar[ll]^{\text{open}}
}
\]
Since $\Gamma$ is strongly convex, $\Gamma^\vee(R)$ is not contained in any hyperplane,
so the monoid $\Gamma^\vee(R)$ generates $\mathfrak{X}^*(R)$ as a group.   This implies 
that $R$ is isomorphic to an open subset of the closed subvariety $R_{\Gamma^\vee(R)} =
\spec k[\Gamma^\vee(R)]$ in $X_\Gamma$, so that $f(R) \cong R$ for every torus $R$ of 
$G$.  Thus $f(\bigcup_T T) = \bigcup_T f(T) \cong \bigcup_T T \subset X_\Gamma$.  Yet 
$\bigcup_T T$ contains a dense open subset $\mathcal{O}$ of $G$ (\cite{Springer}, 
Theorem 6.4.5 and Corollary 7.6.4) and $f|_{\bigcup_T T}$ is an isomorphism, so 
$\mathcal{O} \cong f(\mathcal{O})$ is open in $X_\Gamma$.  Thus $f : G \to X_\Gamma$ 
is a birational morphism~\cite{Hart}.  As $f(G)$ is an open orbit in $X_\Gamma$ 
containing an open subset birationally equivalent to $G$, $f(G)$ is isomorphic to $G$ 
and $f : G \to X_\Gamma$ is an affine $G$-embedding.

Finally, consider $\Gamma(X_\Gamma,f(e)) = \{ \gamma \in \mathfrak{X}_*(G) : 
\limit \gamma(t) f(e) \text{ exists in } X_\Gamma \} = \{ \gamma \in 
\mathfrak{X}_*(G) : \gamma^\circ(A_\Gamma) \subset k[t] \}$, that is, the dual 
map $\gamma^\circ : A_\Gamma \subset k[G] \to k[t,t^{-1}]$ factors through $k[t]$.
Thus, if $\gamma \in \Gamma$, so $v_\gamma(f) \geq 0$ for all $f \in A_\Gamma$, then
$\gamma^\circ(A_\Gamma) \subset k[t]$, and hence $\Gamma \subset \Gamma(X_\Gamma,f(e))$.
Now let $\gamma \in \Gamma(X_\Gamma,f(e))$.  Suppose $\gamma$ is a one-parameter
subgroup of a torus $T$ of $G$.  Then $\limit \gamma(t) f(e)$ exists in the
toric variety $\overline{T} \subset X_\Gamma$, so the classification of toric varieties 
implies that $\gamma \in \Gamma \cap \mathfrak{X}_*(T)$, and hence that 
$\Gamma(X_\Gamma,f(e)) \subset \Gamma$.  Therefore, $X_\Gamma = \spec A_\Gamma$
is a normal affine $G$-embedding such that $\Gamma(X_\Gamma,f(e)) = \Gamma$, which
completes our proof.
\end{proof}



In the next sections, we discuss how our classification in Theorem~\ref{T:classification}
also describes equivariant morphisms between affine $G$-embeddings.  In particular,
we show that equivariant maps between affine $G$-embeddings correspond to inclusions
of the associated strongly convex lattice cones and conversely in 
Section~\ref{SS:morphisms}.  After that, we characterize biequivariant affine 
$G$-embeddings in terms of their cones and, using this result, construct the 
\emph{biequivariant resolution} of an affine $G$-embedding in Section~\ref{SS:biequivariance}.
We remark here that Proposition~\ref{P:biequivariance} below may be seen as an affine 
version of Brion's classification~\cite{B98} of regular $G$-compactifications.

\section{Equivariant morphisms between affine $G$-embeddings}\label{SS:morphisms}

Suppose $f : X \to Y$ is an equivariant morphism between affine $G$-embeddings $X$ 
and $Y$.  By equivariance, if $x_0 \in X$ is a base point for $X$, then $y_0 = 
f(x_0)$ is a base point for $Y$.  Moreover, if $\gamma$ is a one-parameter subgroup 
of $G$ such that $\limit \gamma(t) x_0 = x_\gamma$ exists in $X$, then
$\limit \gamma(t) y_0$ exists in $Y$ and is equal to $f(x_\gamma)$ since 
$f$ is continuous and $f(\gamma(t) x_0) = \gamma(t) f(x_0) = \gamma(t) y_0$ for 
all $t \neq 0$.  Therefore, there is an inclusion $\Gamma(X,x_0) \subset 
\Gamma(Y,f(x_0))$ whenever there exists an equivariant morphism $f : X \to Y$ of 
affine $G$-embeddings.  Moreover, $f$ is the morphism dual to the inclusion of the
subalgebras $A_{\Gamma(X,x_0)} \subset A_{\Gamma(Y,f(x_0))}$ in $k[G]$, because
$f(g \cdot x_0) = g \cdot f(x_0)$ and $Gx_0 = \Omega_X$ is open in $X$ implies
that $f$ is uniquely determined by its value at $x_0$.

\begin{remark}\label{R:nonaffinemorphism}
Suppose $X$ and $Y$ are $G$-embeddings, not necessarily affine.  If $f : X \to Y$ 
is an equivariant morphism and if $x_0$ is a base point for $X$ (i.e., the orbit
$Gx_0$ is isomorphic to $G$ as $G$-varieties), then $f(x_0)$ will be a base point 
for $Y$ and we can define the sets $\Gamma(X,x_0)$ and $\Gamma(Y,f(x_0))$ in the 
same manner as for affine $G$-embeddings.  By the same argument as above, the 
existence of the equivariant morphism $f : X \to Y$ implies that $\Gamma(X,x_0)
\subset \Gamma(Y,f(x_0))$.
\end{remark}

Conversely, suppose $\Gamma(X,x_0) \subset \Gamma(Y,y_0)$ for two affine $G$-embeddings
$X$ and $Y$.  By Theorem~\ref{T:A_X}, $X = \spec A_{\Gamma(X,x_0)}$ and $Y = 
\spec A_{\Gamma(Y,y_0)}$.  The definition of $A_\Gamma = \{ f \in k[G] : v_\gamma(f) 
\geq 0 \text{ for all } \gamma \in \Gamma \}$ implies that $A_{\Gamma(Y,y_0)}
\subseteq A_{\Gamma(X,x_0)}$, so there is a corresponding equivariant morphism of 
$G$-embeddings $X \to Y$ sending $x_0 \mapsto y_0$.  However, by 
Corollary~\ref{C:right-translation}, the subalgebra of $k[G]$ isomorphic to $k[X_\Gamma]$ 
is only determined up to right translations, which correspond to conjugates of the 
cone $\Gamma$.  Thus, 

\begin{prop}\label{P:morphisms}
Suppose $X_1,X_2$ are affine $G$-embeddings and $\Gamma_1,\Gamma_2$ are strongly
convex lattice cones.
\begin{enumerate}
  \item  [a.] If $x_1 \in X_1$ is a base point and $f : X_1 \to X_2$ is an equivariant
         morphism of $G$-embeddings, then $\Gamma(X_1,x_1) \subset \Gamma(X_2,f(x_1))$
         and $f$ is the morphism recovered from the inclusion:
         \[
         \xymatrix{
         k[X_1] \ar[r]^\cong_{\psi_{x_1}^\circ} 
         &
         A_{\Gamma(X_1,x_1)}
         \\
         k[X_2] \ar[r]^\cong_{\psi_{f(x_1)}^\circ} \ar[u]^{f^\circ}
         &
         A_{\Gamma(X_2,x_2)} \ar[u]_\subset
         }
         \]
  \item  [b.] If there is an element $h \in G$ such that $\Gamma_1 \subset h \Gamma_2 h^{-1}$, 
         then there is an equivariant morphism $X_{\Gamma_1} \to X_{\Gamma_2}$ sending 
         the base point $x_1 \in X_{\Gamma_1}$ to $h \cdot x_2 \in X_{\Gamma_2}$, where 
         $\mathfrak{m}_{x_i} = A_{\Gamma_i} \cap \mathfrak{m}_e$ for $i=1,2$. 
\end{enumerate}
\end{prop}

\begin{proof}
We have already proven part 1 of this Lemma in the discussion prior to 
Remark~\ref{R:nonaffinemorphism}.

To prove part 2, by Theorem~\ref{T:classification}, we know that $\Gamma_1,\Gamma_2$ 
correspond to affine $G$-embeddings $X_1 = \spec A_{\Gamma_1}$ and $X_2 = 
\spec A_{\Gamma_2}$, respectively, such that $\Gamma_i = \Gamma(X_i,x_i)$ for $i = 1,2$, 
where $x_i$ is the base point of $X_i$ corresponding to the maximal ideal $\mathfrak{m}_e 
\cap A_{\Gamma_i}$.  Now suppose there is an element $h \in G$ such that $\Gamma_1 \subset 
h \Gamma_2 h^{-1}$ as in the statement of the lemma.  We know that $h \Gamma_2 h^{-1} = 
h \Gamma(X_2,x_2) h^{-1} = \Gamma(X,h \cdot x_2)$ by formula (\ref{E:Gamma(X,hx)}).  
Hence we have $\Gamma(X_1,x_1) \subset \Gamma(X_2,h \cdot x_2)$, so that 
$A_{\Gamma(X_1,x_1)} \supset A_{\Gamma(X_2,h \cdot x_2)}$ and thus $X_1 \to X_2$ exists, 
is equivariant, and maps $x_1 \mapsto h \cdot x_2$ as claimed.
\end{proof}

\begin{eg}[The subcone of $\Gamma(X,x_0)$ associated to a torus 
           closure]\label{EG:Gamma(GbarT)inGamma(X)}
If $X$ is an affine $G$-embedding and $T$ is a maximal torus of $G$ whose closure 
in $X$ is denoted $\overline{T}$, then $G \times^T \overline{T}$ is an affine 
$G$-embedding and we have an equivariant morphism $G \times^T \overline{T} \to X$.  
Let $\Gamma = \Gamma(X,x_0)$ and $\sigma = \Gamma(X,x_0) \cap \mathfrak{X}_*(T)$.  
Then 
\[
\Gamma(G \times^T \overline{T},[e,x_0]) 
= \bigcup_{P \supset T} P \bullet (\sigma \cap \Delta_P(T)) ,
\] 
where the union is taken over all parabolic subgroups of $G$ containing $T$ and $\Delta_P(T) 
= \{ \gamma \in \mathfrak{X}_*(T) : P(\gamma) \supseteq P \}$.  There are only finitely many parabolic subgroups $P$ containing $T$, as there are only a finite number of Borel subgroups containing $T$ (in one-to-one correspondence with the elements of the Weyl group $W(T,G)$, [\cite{Springer}, Corollary 6.4.12]), and, for each Borel subgroup $B$, there are only finitely many parabolic subgroups $P \supset B$ (indexed by subsets of the basis of positive roots of $T$ relative to $B$, [\cite{Springer}, Theorem 8.4.3]).  In contrast, $\Gamma = \bigcup_{Q \subset G} Q \bullet (\Gamma \cap \Delta_Q(G))$, where the union is taken over the set of all parabolic subgroups $Q$ of $G$ (there are infinitely many, as there are infinitely many Borel subgroups each corresponding to cosets in $G/B$, which is projective) and $\Delta_Q(G) = \{ \gamma \in \mathfrak{X}_*(G) : P(\gamma) \supseteq Q \}$.  Clearly, for each $P \supset T$, 
\[
\sigma \cap \Delta_P(T) 
= (\Gamma \cap \mathfrak{X}_*(T)) \cap \Delta_P(T) 
= \Gamma \cap \Delta_P(T) 
= \Gamma \cap \Delta_P(G)
\]
since $\Delta_P(G) \subset \mathfrak{X}_*(T')$ for all maximal tori $T'$ contained in 
$P$.  Therefore, $\Gamma(G \times^T \overline{T},[e,x_0])$ is a finite union of the 
``parabolic components'' of $\Gamma(X,x_0)$, so $\Gamma(G \times^T \overline{T},[e,x_0])$ 
is a finite polysimplicial subcomplex of $\Gamma(X,x_0)$, as the $\Delta_Q(G)$ provide 
a triangulation of $\mathfrak{X}_*(G)$~\cite{GIT}.  
\end{eg}


\section{Biequivariant resolutions}\label{SS:biequivariance}

In this section, we show that every affine $G$-embedding $X$ canonically determines
a $(G \times G)$-equivariant affine $G$-embedding $X_G$ together with a 
left-$G$-equivariant morphism $X_G \to X$, which we call the biequivariant resolution
of $X$.  As we will be working with varieties some of which only have a left action
and others both a left and a right action, we will be careful to specify how $G$
acts on varieties discussed in this section.

Our first tool is the following proposition, which allows us to detect when an affine
$G$-embedding also has a right-$G$-action $X \times G \to X$ extending the 
multiplication in $G$.

\begin{prop}\label{P:biequivariance}
An affine $G$-embedding $X$ will have both a left and a right $G$-action, and thus be 
a biequivariant $G$-embedding, if and only if the associated strongly convex lattice 
cone $\Gamma(X,x_0)$, for any choice of base point $x_0 \in X$, is $G$-stable for the 
conjugation action of $G$ on $\mathfrak{X}_*(G)$.
\end{prop}

\begin{proof}
Suppose that $X$ is a $(G \times G)$-equivariant affine $G$-embedding and let $x \in X$ 
be a base point.  Then $X$ is determined by the strongly convex lattice cone $\Gamma(X,x) 
= \{ \gamma \in \mathfrak{X}_*(G) : \limit \gamma(t) x \text{ exists in } X \}$.  
Let $h \in G$ and recall that $h \Gamma(X,x) h^{-1} = \Gamma(X,h \cdot x)$ by formula 
(\ref{E:Gamma(X,hx)}).  Thus it is enough to show that $\gamma \in \Gamma$ if and only 
if $\gamma \in \Gamma(X,h \cdot x)$.  Assume $\limit \gamma(t) x$ exists in $X$.  
Then consider $\limit [ \gamma(t) \cdot hx ] = \limit [ \gamma(t) \cdot xh' ] 
= [ \limit \gamma(t) \cdot x ] \cdot h'$, for some $h' \in G$, and this limit 
exists in $X$ by Remark~\ref{R:bi-limits}, since there is a right $G$-action on $X$.  
Thus $\Gamma(X,x) \subset h \Gamma(X,x) h^{-1}$.  Now assume that $\gamma' \in 
\Gamma(X,h \cdot x)$.  Then, the same argument implies $\gamma' \in \Gamma(X,x)$ using 
Remark~\ref{R:bi-limits}.  Therefore, for every $h \in G$, $\Gamma(X,x) = h \Gamma(X,x) 
h^{-1}$.  Thus $\Gamma(X,x)$ is $G$-stable for the conjugation action of $G$ on 
$\mathfrak{X}_*(G)$ for any choice of base point $x \in X$.

Now assume that $\Gamma$ is a strongly convex lattice cone which is $G$-stable for 
the conjugation action of $G$ on $\mathfrak{X}_*(G)$.  Theorem~\ref{T:classification} 
implies that $\Gamma = \Gamma(X_\Gamma,x_0)$ for some base point $x_0 \in X_\Gamma = 
\spec A_\Gamma$.  Then, by Corollary~\ref{C:right-translation}, for any $h \in G$ we 
have $r_h(A_{\Gamma(X_\Gamma,x_0)}) = A_{\Gamma(X_\Gamma,h \cdot x_0)} = 
A_{h \Gamma(X_\Gamma,x_0) h^{-1}} = A_\Gamma$.  Hence $A_\Gamma$ is a left- and
right-$G$-invariant subalgebra of $k[G]$, so there is a right $G$-action on 
$X_\Gamma = \spec A_\Gamma$, which extends the multiplication of $G$ by~\cite{Springer}, 
Proposition 2.3.6.  Thus $X_\Gamma$ is a $(G \times G)$-equivariant affine $G$-embedding.
\end{proof}

\begin{coro}\label{C:biequivariance}
If $X$ is a biequivariant affine $G$-embedding and $T$ is any maximal torus of $G$,
then the closure of $T$ in $X$ determines $X$ completely.
\end{coro}

\begin{proof}
Let $X$ be a biequivariant affine $G$-embedding with lattice cone $\Gamma = \Gamma(X,x_0)$ 
for some choice of base point in $X$. Let $T$ be any maximal torus of $G$.  Consider 
$\overline{T} \subset X$, which is an affine toric variety for $T$.  Let $\sigma \subset 
\mathfrak{X}_*(T)$ be the cone of one-parameter subgroups of $T$ with specializations in 
$\overline{T}$.  Then $\sigma = \Gamma \cap \mathfrak{X}_*(T)$, by definition of $\Gamma$.  
If $T'$ is any other maximal torus of $G$, then $T' = gTg^{-1}$ for some $g \in G$ and 
$\Gamma \cap \mathfrak{X}_*(T') = g[g^{-1}(\Gamma \cap \mathfrak{X}_*(T'))g]g^{-1} = 
g[ g^{-1} \Gamma g \cap g^{-1} \mathfrak{X}_*(T') g ] g^{-1} = g[ \Gamma \cap 
\mathfrak{X}_*(T) ]g^{-1} = g \sigma g^{-1}$.  Thus $\Gamma = \bigcup_{g \in G} g \sigma 
g^{-1}$ is determined by the cone $\sigma$, so $X$ is determined by its toric
subvariety $\overline{T} = T_\sigma$.
\end{proof}

We note the similarity between our classification of biequivariant affine 
$G$-embeddings and Brion's classification of regular $G$-compactifications.  
In~\cite{B98}, Brion classified regular compactifications of a group $G$ by the 
$W(T,G)$-invariant fan associated to the closure of any maximal torus $T$ of $G$ 
in the compactification, demonstrating that regular compactifications of $G$ are 
completely determined by any of the associated toric subvarieties. 
Likewise, Corollary~\ref{C:biequivariance} implies that if an affine $G$-embedding $X$ 
is $(G \times G)$-equivariant and $T$ is a maximal torus of $G$, then the cone for
$\overline{T}$ recovers $X$.  Therefore, we may think of Proposition~\ref{P:biequivariance} 
and Corollary~\ref{C:biequivariance} as an affine version of Brion's classification.  

\begin{remark}\label{R:basepointfree}
Suppose $X$ is a $(G \times G)$-equivariant affine $G$-embedding and suppose that $x_0 \in 
X$ is a base point.  By Proposition~\ref{P:biequivariance} above, the set $\Gamma(X,x_0)$
is stable under the conjugation action of $G$.  However, by formula (\ref{E:Gamma(X,hx)}),
we conclude that $\Gamma(X,x_0) = h\Gamma(X,x_0)h^{-1} = \Gamma(X,h \cdot x_0)$ for all 
$h \in G$.  Therefore, the strongly convex lattice cone associated to a biequivariant
affine $G$-embedding is independent of the choice of base point.  
\end{remark}

We use these observations to construct the biequivariant resolution of
an affine $G$-embedding $X$.  Suppose $X$ is an affine $G$-embedding,
$x_0 \in X$ is a base point, and $\Gamma = \Gamma(X,x_0)$ is the
corresponding strongly convex lattice cone.  Then there is a unique
maximal $G$-stable subset
\begin{equation}\label{E:Gamma^G}
\Gamma^G := \bigcap_{h \in G} h \Gamma h^{-1}
\end{equation}
of $\Gamma$.  In fact,

\begin{lemma}\label{L:Gamma^G}
If $\Gamma$ is a strongly convex lattice cone in $\mathfrak{X}_*(G)$,
then $\Gamma^G$ is a strongly convex lattice cone.
\end{lemma}

\begin{proof}
Let $\Gamma$ be a strongly convex lattice cone.  Let $\Xi$ be a Kempf
state such that $\Gamma = \Xi^\vee$.  Define $\Gamma^G = \bigcap_{h
\in G} h \Gamma h^{-1}$ as in (\ref{E:Gamma^G}).  We claim that
$G_*\Xi : R \mapsto \bigcup_{g \in G} g_*\Xi(R)$ is a Kempf state
and that $\Gamma^G = (G_*\Xi)^\vee$.

First, $G_*\Xi$ is a state, for each $g_*\Xi$ is a state implies that 
$G_*\Xi(R)$ is nonempty for all $R$ and that if $S \supset R$, then $G_*\Xi(R) = 
\bigcup_{g \in G} g_*\Xi(R) = \bigcup_{g \in G} \Res^S_R[g_*\Xi(S)] = 
\Res^S_R\left[\bigcup_{g \in G} g_*\Xi(S) \right] = \Res^S_R[ G_*\Xi(S) ]$.  Next, 
$G_*\Xi$ is bounded, for if $h \in G$, then $h_*[G_*\Xi](R) = \bigcup_{g \in G} 
h_*[g_*\Xi(R)] = \bigcup_{g \in G} (hg)_*\Xi(R) = G_*\Xi(R)$, so $\bigcup_{h \in G} 
h_*[G_*\Xi](R) = G_*\Xi(R) = \bigcup_{g \in G} g_*\Xi(R)$ is a finite subset of 
$\mathfrak{X}^*(R)$ for each $R$, since $\Xi$ is a bounded state.
Consider the numerical function $\mu(G_*\Xi)$.  If $\gamma \in
\mathfrak{X}_*(G)$, then $\mu(G_*\Xi,\gamma) = \min_{\chi \in
G_*\Xi(\gamma)} \langle \chi,\gamma \rangle = \min_{\chi \in \bigcup
g_*\Xi(\gamma)} \langle \chi,\gamma \rangle = \min_{g \in G}
\min_{\chi \in g_*\Xi(\gamma)} \langle \chi,\gamma \rangle = 
\min_{g \in G} \mu(g_*\Xi,\gamma)$.  Therefore, 
\begin{equation}\label{E:muG*Xi}
\mu(G_*\Xi,\gamma) = \min_{g \in G} \mu(g_*\Xi,\gamma).
\end{equation}
Each $g_*\Xi$ is a Kempf state, hence admissible.  Hence, if $\gamma$
is a one-parameter subgroup of $G$ and $p \in P(\gamma)$ belongs to
the parabolic subgroup associated to $\gamma$, then $\mu(G_*\Xi,p
\bullet \gamma) = \min_{g \in G} \mu(g_*\Xi,p \bullet \gamma) =
\min_{g \in G} \mu(g_*\Xi,\gamma) = \mu(G_*\Xi,\gamma)$.  Thus
$G_*\Xi$ is admissible, so $G_*\Xi$ is a Kempf state.  Moreover,
(\ref{E:muG*Xi}) implies $\gamma \in (G_*\Xi)^\vee$ if and only if
$\mu(g_*\Xi,\gamma) \geq 0$ for all $g \in G$.  But $\{ \gamma \in
\mathfrak{X}_*(G) : \mu(g_*\Xi,\gamma) \geq 0 \} = (g_*\Xi)^\vee = g
(\Xi^\vee) g^{-1} = g \Gamma g^{-1}$ by (\ref{E:g*Xi}), so
$(G_*\Xi)^\vee = \Gamma^G$ as claimed.

Thus $\Gamma^G = \{ \gamma \in \mathfrak{X}_*(G) : \mu(G_*\Xi,\gamma)
\geq 0 \}$, which implies $\gamma \in \Gamma^G$ if and only if
$\gamma^n \in \Gamma^G$ for every positive integer $n > 0$, as
$\mu(G_*\Xi,\gamma^n) = n \mu(G_*\Xi,\gamma)$.  We use this condition
to demonstrate that $\Gamma^G$ is saturated with respect to the
equivalence relation $\sim$.  Suppose $\gamma \in \Gamma^G$ and
$\delta \sim \gamma$.  Then there are positive integers $m,n > 0$ and
an element $g \in P(\gamma)$ such that $\delta^m = g \gamma^n g^{-1}$.
Yet $\gamma \in \Gamma^G$ implies that $\gamma^n \in \Gamma^G$, for $n
> 0$.  Hence $\delta^m = g \gamma^n g^{-1} \in g \Gamma^G g^{-1} = g
\left(\bigcap_{h \in G} h \Gamma h^{-1} \right) g^{-1} = \Gamma^G$.
Since $m > 0$, this implies that $\delta \in \Gamma^G$ as well, so
that $\Gamma^G$ is saturated.

To prove that $\Gamma^G$ is a convex lattice cone, it is left to show
that the set $\Gamma^G_1/\sim$ is finite.  
This is clear since $\Gamma^G_1 \subset \Gamma_1$ and $\Gamma$ is a
convex lattice cone.

Finally, as $\Gamma^G \subset \Gamma$, it is clear that
$\gamma,\gamma^{-1} \in \Gamma^G$ implies that $\gamma = \varepsilon$
is the trivial one-parameter subgroup of $G$.  Furthermore,
$\varepsilon \in h \Gamma h^{-1}$ for all $h \in G$, so $\varepsilon
\in \Gamma^G$.  Thus $\gamma, \gamma^{-1} \in \Gamma^G$ if and only if
$\gamma = \varepsilon$.  Therefore, $\Gamma^G$ is a strongly convex
lattice cone in $\mathfrak{X}_*(G)$.
\end{proof}

Let $X$ be an affine $G$-embedding and select a base point $x_0 \in X$.  Let 
$\Gamma = \Gamma(X,x_0)$ and define $\Gamma^G = \bigcap_{h \in G} h \Gamma h^{-1}$,
which is a strongly convex lattice cone.  Then $X_G := \spec A_{\Gamma^G}$
is a $(G \times G)$-equivariant affine $G$-embedding by Theorem~\ref{T:classification}
and Proposition~\ref{P:biequivariance}, and there is a left-$G$-equivariant morphism 
$\beta_{(X,x_0)} : X_G \to X$ corresponding to the inclusion $A_{\Gamma(X,x_0)} 
\subset A_{\Gamma^G}$ of subalgebras in $k[G]$.  

\begin{definition}\label{D:biresolution}
We call $X_G$ together with the left-$G$-equivariant morphism
$\beta_{(X,x_0)} : X_G \to X$ the \emph{biequivariant resolution} of $X$.
\end{definition}

We make a few remarks about the biequivariant resolution of an affine $G$-embedding
before we prove its universal property.  First, if $X$ is already a $(G \times
G)$-equivariant affine $G$-embedding, then $\Gamma = \Gamma^G$, so $X_G = X$.
However, it is possible in some cases that $X_G = G$ even when $X \neq G$.

Second, while the cone $\Gamma^G$ is canonical, the morphism $\beta_{(X,x_0)} : 
X_G \to X$ is only defined up to right translation, which corresponds to a different 
choice of base point $x \in X$ as follows: $r_h(\beta_{(X,x_0)}) = \beta_{(X,h \cdot 
x_0)} : X_G \to X$.  This is true because the identification of $k[X]$ with a subalgebra 
of $k[G]$ is determined by the selection of a base point and changing the base point 
corresponds to right translation of the subalgebra in $k[G]$ by 
Corollary~\ref{C:right-translation}.  Then the morphism $\beta_{(X,h \cdot x_0)} : X_G 
\to X$ is defined by the inclusion of $A_{\Gamma(X,h \cdot x_0)} = r_h(A_{\Gamma(X,x_0)}) 
\subset A_{\Gamma^G}$, from which it is clear that $\beta_{(X,h \cdot x_0)} = 
r_h(\beta_{(X,x_0)})$.


\begin{thm}[Universal Property of Biequivariant Resolutions]\label{T:universalbi}
Suppose $X$ is a left-$G$-equivariant affine $G$-embedding, $x_0 \in X$ is a base point,
$Y$ is any $(G \times G)$-equivariant affine $G$-embedding, and $\varphi : Y \to X$
is a left-$G$-equivariant morphism.  Then there is a unique $(G \times G)$-equivariant 
morphism $\varphi_{(X_G,x_0)} : Y \to X_G$ such that $\varphi = \beta_{(X,x_0)} \circ 
\varphi_{(X_G,x_0)}$.  That is, the biequivariant resolution of $X$ satisfies the 
following diagram:
\[
\xymatrix{
Y \ar[rr]^{\forall \, \varphi} \ar@{-->}[dr]_{\exists! \, \varphi_{(X_G,x_0)}}
&&
X
\\
&
X_G \ar[ur]_{\beta_{(X,x_0)}}
}
\]
If $x_1 = h \cdot x_0$ is another base point for $X$, then $\beta_{(X,x_1)} = 
r_h(\beta_{(X,x_0)})$ and $\varphi_{(X_G,x_1)} = r_{h^{-1}}(\varphi_{(X_G,x_0)})$.
\end{thm}

\begin{proof}
Let $X$ be an affine $G$-embedding and identify $k[X]$ with the left-invariant subalgebra 
$A_{\Gamma(X,x_0)}$ of $k[G]$ by selecting a base point $x_0 \in X$.  Suppose that $Y$ 
is a $(G \times G)$-equivariant affine $G$-embedding and that $\varphi : Y \to X$ is a 
left-$G$-equivariant morphism from $Y$ to $X$.  By equivariance, there is a unique 
element $y_0 \in Y$ such that $\varphi(y_0) = x_0$, since $\varphi|_{\Omega_Y} : 
\Omega_Y \to \Omega_X$ is an isomorphism with $\Omega_Y \cong \Omega_X \cong G$.  Then 
the cone $\Gamma(Y,y_0)$ is $G$-stable by Proposition~\ref{P:biequivariance} and is a 
subset of $\Gamma(X,x_0)$ by Proposition~\ref{P:morphisms}.  Using $y_0$, identify 
$k[Y]$ with the $(G \times G)$-invariant subalgebra $A_{\Gamma(Y,y_0)}$ of $k[G]$.  By
Corollary~\ref{C:right-translation}, $A_{\Gamma(Y,h \cdot y_0)} = r_h(A_{\Gamma(Y,y_0)}) 
= A_{\Gamma(Y,y_0)}$, for all $h \in G$, so the algebra $A_{\Gamma(Y,y_0)}$ is 
independent of the choice of base point.  Hence it is the only subalgebra of $k[G]$ 
which is isomorphic to $k[Y]$ as a $(G \times G)$-algebra by Theorem~\ref{T:classification}.

Clearly $\Gamma(Y,y_0) = \bigcap_{h \in G} h \Gamma(Y,y_0) h^{-1} \subset 
\bigcap_{h \in G} h \Gamma(X,x_0) h^{-1} = \Gamma(X,x_0)^G \subset \Gamma(X,x_0)$, so 
that $A_{\Gamma(X,x_0)} \subset A_{\Gamma(X,x_0)^G} \subset A_{\Gamma(Y,y_0)}$, as
indicated in the diagram below.  Moreover, as $\Gamma(X,x_0)^G$ is $G$-stable, 
applying Corollary~\ref{C:right-translation} again shows that $A_{\Gamma(X,x_0)^G}$
is a uniquely determined subalgebra of $k[G]$ which is independent of the choice of
base point $x_0 \in X$ made above.  Thus we have the following diagram of algebras:
\[
\xymatrix{
&
k[G]
&
\\
k[Y] \ar[r]^\cong_{\psi_{y_0}^\circ}
&
A_{\Gamma(Y,y_0)} \ar[u]_\subset
&
\\
&
&
A_{\Gamma(X,x_0)^G} = k[X_G] \ar[ul]_\subset
\\
k[X] \ar[uu]^{\varphi^\circ} \ar[r]^\cong_{\psi_{x_0}^\circ}
&
A_{\Gamma(X,x_0)} \ar[uu]^\subset \ar[ur]_\subset
}
\]
The morphism $\varphi_{(X_G,x_0)} : Y \to X_G$ corresponds to the restriction of the 
homomorphism $(\psi_{y_0}^\circ)^{-1}$ from $k[X_G] = A_{\Gamma(X,x_0)^G} \to k[Y]$, 
while $\beta_{(X,x_0)} : X_G \to X$ is dual to the composition of $\psi_{x_0}^\circ : 
k[X] \to A_{\Gamma(X,x_0)}$ followed by the canonical inclusion $A_{\Gamma(X,x_0)} 
\subset A_{\Gamma(X,x_0)^G}$.  Then it is clear that $\varphi$ factors as $\varphi = 
\beta_{(X,x_0)} \circ \varphi_{(X_G,x_0)}$.

Now suppose that $x_1 = h \cdot x_0$ is another base point in $X$.  Then $h \cdot y_0$
is the unique element of $Y$ such that $\varphi(y) = x_1$, for $\varphi(h \cdot y_0)
= h \cdot \varphi(y_0) = h \cdot x_0$.  Then $\psi_{x_1}^\circ$ is an isomorphism
from $k[X]$ to $A_{\Gamma(X,h \cdot x_0)} = A_{h \Gamma(X,x_0) h^{-1}} = 
r_h(A_{\Gamma(X,x_0)}) \subset k[G]$.  Yet $\Gamma(Y,h \cdot y_0) = \Gamma(Y,y_0)
\subset \Gamma(X,x_0)^G = \Gamma(X,h \cdot x_0)^G \subset \Gamma(X,h \cdot x_0)$,
so we have $\varphi^\circ = (\psi_{h \cdot y_0}^\circ)^{-1} \circ \psi_{h \cdot x_0}^\circ
: k[X] \to k[Y]$.  As above, this factors through the inclusions $A_{\Gamma(X,h \cdot
x_0)} \subset k[X_G] = A_{\Gamma(X,h \cdot x_0)^G} \subset A_{\Gamma(Y,h \cdot y_0)}$, 
giving morphisms $\varphi_{(X_G,h \cdot x_0)} : Y \to X_G$ and $\beta_{(X,h \cdot x_0)}
: X_G \to X$ corresponding to $(\psi_{h \cdot y_0}^\circ)^{-1} = 
[r_h(\psi_{y_0}^\circ)]^{-1} = r_{h^{-1}}[(\psi_{y_0}^\circ)^{-1}]$ and 
$\psi_{h \cdot x_0}^\circ = r_h(\psi_{x_0}^\circ)$, respectively.
Hence $\varphi_{(X_G,h \cdot x_0)} = r_{h^{-1}}(\varphi_{(X_G,x_0)})$ and 
$\beta_{(X,h \cdot x_0)} = r_h(\beta_{(X,x_0)})$.
\end{proof}


\end{document}